\input amstex\documentstyle{amsppt}  
\pagewidth{12.5cm}\pageheight{19cm}\magnification\magstep1  
\topmatter
\title The flag manifold over the semifield $\bold Z$\endtitle
\author G. Lusztig\endauthor
\address{Department of Mathematics, M.I.T., Cambridge, MA 02139}\endaddress
\thanks{Supported by NSF grant DMS-1855773.}\endthanks
\endtopmatter   
\document
\define\rhd{\rightharpoondown}

\define\dw{\dot w}

\define\ds{\dot s}

\define\mpb{\medpagebreak}

\define\frl{\forall}

\define\si{\sim}
\define\wt{\widetilde}
\define\sqc{\sqcup}

\define\qua{\quad}

\define\part{\partial}
\define\emp{\emptyset}

\define\ra{\rangle}

\define\m{\mapsto}
\define\do{\dots}
\define\la{\langle}

\define\sub{\subset}    

\define\T{\times}
\define\ti{\tilde}
\define\nl{\newline}
\redefine\i{^{-1}}

\define\un{\underline}

\define\ot{\otimes}

\define\End{\text{\rm End}}

\define\supp{\text{\rm supp}}

\redefine\b{\beta}
\redefine\c{\chi}
\define\g{\gamma}
\redefine\d{\delta}
\define\e{\epsilon}

\redefine\o{\omega}
\define\p{\pi}
\define\ph{\phi}

\define\s{\sigma}
\redefine\t{\tau}
\define\th{\theta}
\define\k{\kappa}
\redefine\l{\lambda}
\define\z{\zeta}
\define\x{\xi}

\redefine\G{\Gamma}

\define\Om{\Omega}

\define\ii{\bold i}

\define\kk{\bold k}

\define\CC{\bold C}

\define\NN{\bold N}

\define\RR{\bold R}

\define\ZZ{\bold Z}

\define\cb{\Cal B}

\define\ce{\Cal E}

\define\ci{\Cal I}
\define\cj{\Cal J}

\define\cl{\Cal L}

\define\co{\Cal O}

\define\cs{\Cal S}
\define\ct{\Cal T}

\define\cz{\Cal Z}
\define\cx{\Cal X}
\define\cy{\Cal Y}

\define\fG{\frak G}

\define\ta{\ti a}

\define\tA{\ti A}

\define\tF{\ti F}
\define\tG{\ti G}
\define\tH{\ti H}

\define\tP{\ti P}

\define\tV{\ti V}

\define\tX{\ti X}

\define\bE{\bar E}

\define\bul{\bullet}

\define\cir{\circ}

\define\tcb{\ti{\cb}}

\head Introduction\endhead
\subhead 0.1\endsubhead
Let $G$ be a connected semisimple simply connected algebraic group over $\CC$ with a fixed pinning (as in
\cite{L94b, 1.1}). In this paper we assume that $G$ is of simply laced type.
Let $\cb$ be the variety of Borel subgroups of $G$. In \cite{L94b, 2.2, 8.8} a submonoid $G_{\ge0}$
of $G$ and a subset $\cb_{\ge0}$ of $\cb$ with an action of $G_{\ge0}$ (see \cite{L94b, 8.12})
was defined. (When $G=SL_n$, $G_{\ge0}$ is the submonoid consisting of the
real, totally positive matrices in $G$.) More generally, for any semifield $K$,
a monoid $\fG(K)$ was defined in \cite{L19a}, so that when $K=\RR_{>0}$ we have $\fG(K)=G_{\ge0}$.
(In the case where $K$ is $\RR_{>0}$ or the semifield in (i) or (ii) below, a monoid $G(K)$ already 
appeared in \cite{L94b, 2.2, 9.10}; it was identified with $\fG(K)$ in \cite{L19b}.) 

This paper is concerned with the question of definining the flag manifold over a semifield $K$
with an action of the monoid $\fG(K)$ so that in the case where $K=\RR_{>0}$ we recover $\cb_{\ge0}$
with its $G_{\ge0}$-action.

In \cite{L19b, 4.9}, for any semifield $K$, a definition of the flag manifold over $K$
was given (based on ideas of Marsh and Rietsch \cite{MR}); but in that 
definition the lower and upper triangular part of $G$ play an asymmetric role
 and as a consequence only a part of $\fG(K)$ acts on $\cb(K)$ (unlike the case
$K=\RR_{>0}$ when the entire $\fG(K)$ acts). To get the entire $\fG(K)$ act one needs a conjecture
stated in \cite{L19b, 4.9} which is still open.

In this paper we get around that conjecture and provide an unconditional definition of the flag
manifold  (denoted by $\cb(K)$) 
over a semifield $K$ with an action of $\fG(K)$ assuming that $K$ is either

(i) the semifield consisting of all rational functions in $\RR(x)$ (with $x$ an indeterminate) of the form
$x^ef_1/f_2$ where $e\in\ZZ$ and $f_1\in\RR[x],f_2\in\RR[x]$ have constant term in $\RR_{>0}$ (standard sum and 
product); or

(ii) the semifield $\ZZ$ in which the sum of $a,b$ is $\min(a,b)$ and the product of $a,b$ is $a+b$.
\nl
For $K$ as in (i) we give two definitions of $\cb(K)$; one of them is elementary and the other is less so,
being based on the theory of canonical bases (the two definitions are shown to be equivalent). 
For $K$ as in (ii) we only give a definition based on the theory of canonical bases.

A part of our argument involves a construction of an analogue of the finite dimensional irreducible
representations of $G$ when $G$ is replaced by the monoid $\fG(K)$ where $K$ is any semifield.

Let $W$ be the Weyl group of $G$. Now $W$ is naturally a Coxeter group with generators
$\{s_i;i\in I\}$ and length function $w\m|w|$. 
Let $\le$ be the Chevalley partial order on $W$.

In \S3 we prove the following result which is a $\ZZ$-analogue of a result (for $\RR_{>0}$) in \cite{MR}.

\proclaim{Theorem 0.2} The set $\cb(\ZZ)$ has a canonical partition into pieces $P_{v,w}(\ZZ)$ indexed by the
pairs $v\le w$ in $W$. Each such piece $P_{v,w}(\ZZ)$ is in bijection with $\ZZ^{|w|-|v|}$; in fact, there is an
explicit  bijection $\ZZ^{|w|-|v|}@>\si>>P_{v,w}(\ZZ)$ for any reduced expression of $w$.
\endproclaim

In \S3 we also prove a part of a conjecture in \cite{L19b, 2.4} which attaches to any $v\le w$ in $W$
a certain subset of a canonical basis, see 3.10.

In \S4 we show that our definitions do not depend on the choice of a (very dominant) weight $\l$.

In \S5 we show how some of our results extend to the non-simply laced case.

\head Contents\endhead

1. Definition of $\cb(\ZZ)$.

2. Preparatory results.

3. Parametrizations.

4. Independence on $\l$.

5. The non-simply laced case.

\head 1. Definition of $\cb(\ZZ)$\endhead
\subhead 1.1\endsubhead
In this section we will give the definition of the flag manifold $\cb(K)$ when $K$ is as in 0.1(i),(ii).

\subhead 1.2\endsubhead
We fix some notation on $G$. Let $w_I$ be the longest element of $W$. For $w\in W$ let 
$\ci_w$ be the set of all sequences $\ii=(i_1,i_2,\do,i_m)$ in $I$ such that $w=s_{i_1}s_{i_2}\do s_{i_m}$, 
$m=|w|$. 

The pinning of $G$ consists of two opposed Borel subgroups $B^+,B^-$ with unipotent radicals $U^+,U^-$ and 
root homomorphisms $x_i:\CC@>>>U^+$, $y_i:\CC@>>>U^-$ indexed by $i\in I$. Let $T=B^+\cap B^-$, a maximal torus.
Let $\cy$ be the group of one parameter subgroups $\CC^*@>>>T$; let 
$\cx$ be the group of characters $T@>>>\CC^*$. Let $\la,\ra:\cy\T\cx@>>>\ZZ$ be the canonical pairing.
The simple coroot corresponding to $i\in I$ is denoted again by $i\in\cy$; let $i'\in\cx$ 
be the corresponding simple root. Let $\cx^+=\{\l\in\cx;\la i,\l\ra\ge0\qua\frl i\in I\}$,
$\cx^{++}=\{\l\in\cx;\la i,\l\ra\ge1\qua\frl i\in I\}$. 
Let $G(\RR)$ be the subgroup of $G$ generated by $x_i(t),y_i(t)$ with $i\in I,t\in\RR$.
Let $\cb(\RR)$ be the subset of $\cb$ consisting of all $B\in\cb$ such that $B=gB^+g\i$ for some 
$g\in G(\RR)$. We have $G_{\ge0}\sub G(\RR)$, $\cb_{\ge0}\sub\cb(\RR)$.
For $i\in I$ we set $\ds_i=y_i(1)x_i(-1)y_i(1)\in G(\RR)$, an element normalizing $T$.
For $(B,B')\in\cb\T\cb$ we write $pos(B,B')$ for the relative position of $B,B'$ (an element of $W$).

\subhead 1.3\endsubhead
Let $K$ be a semifield. Let $K^!=K\sqc\{\circ\}$ where $\circ$ is a symbol.
We extend the sum and product on $K$ to a sum and product on $K^!$ by definining
$\circ+a=a$, $a+\circ=a$, $\circ\T a=\circ, a\T\circ=\circ$ for $a\in K$ and
$\circ+\circ=\circ, \circ\T\circ=\circ$. Thus $K^!$ becomes a monoid under addition and a monoid
under multiplication. Moreover the distributivity law holds on $K^!$.
When $K$ is $\RR_{>0}$ we have $K^!=\RR_{\ge0}$ with $\circ=0$ and the usual sum and product.
When $K$ is as in 0.1(i), $K^!$ can be viewed as the subset of $\RR(x)$ given by $K\cup\{0\}$ with $\circ=0$
and the usual sum and product. When $K$ is as in 0.1(ii) we have $0\in K$ and $\circ\ne0$.

\subhead 1.4\endsubhead
Let $V={}^\l V$ be the finite dimensional simple $G$-module over $\CC$ with highest weight $\l\in\cx^+$. For 
$\nu\in\cx$ 
let $V_\nu$ be the $\nu$-weight space of $V$ with respect to $T$. Thus $V_\l$ is a line. We fix 
$\x^+={}^\l\x^+$ in $V_\l-0$.
For each $i\in I$ there are well defined linear maps $e_i:V@>>>V,f_i:V@>>>V$
such that $x_i(t)\x=\sum_{n\ge0}t^ne_i^{(n)}\x,y_i(t)\x=\sum_{n\ge0}t^nf_i^{(n)}\x$ for $\x\in V,t\in\CC$.
Here $e_i^{(n)}=(n!)\i e_i^n:V@>>>V, f_i^{(n)}=(n!)\i f_i^n:V@>>>V$ are zero for $n\gg0$. For an integer
$n<0$ we set $e_i^{(n)}=0,f_i^{(n)}=0$.

Let $\b={}^\l\b$ be the canonical basis of $V$ (containing $\x^+$) defined in \cite{L90a}. 
Let $\x^-$ be the lowest weight vector in $V-0$ contained in $\b$. 
For $b\in\b$ we have $b\in V_{\nu_b}$ for a well defined $\nu_b\in\cx$, said to be the weight of $b$.
By a known property of $\b$ (see \cite{L90a, 10.11} and \cite{L90b,\S3}, or alternatively \cite{L93, 22.1.7}), 
for $i\in I,b\in\b,n\in\ZZ$ we have 
$$e_i^{(n)}b=\sum_{b'\in\b}c_{b,b',i,n}b',\qua f_i^{(n)}b=\sum_{b'\in\b}d_{b,b',i,n}b'$$
where 
$$c_{b,b',i,n}\in\NN,\qua d_{b,b',i,n}\in\NN.$$
 Hence for $i\in I,b\in\b,t\in\CC$ we have 
$$x_i(t)b=\sum_{b'\in\b,n\in\NN}c_{b,b',i,n}t^nb',\qua y_i(t)b=\sum_{b'\in\b,n\in\NN}d_{b,b',i,n}t^nb'.$$
For any $i\in I$ there is a well defined function $z_i:\b@>>>\ZZ$ such that for $b\in\b$, $t\in\CC^*$ we have 
$i(t)b=t^{z_i(b)}b$.

\mpb

Let $P={}^\l P$ be the variety of $\CC$-lines in $V$. 
Let $P^\bul={}^\l P^\bul$ be the set of all $L\in P$ such that for some $g\in G$ we have $L=gV_\l$. 
Now $P^\bul$ is a closed subvariety of $P$. For any $L\in P^\bul$ let $G_L=\{g\in G;gL=L\}$; this is a
parabolic subgroup of $G$. 

Let $V^\bul={}^\l V^\bul=\cup_{L\in P^\bul}L$, a closed subset of $V$. 
For any $\x\in V,b\in\b$ we define $\x_b\in\CC$  by $\x=\sum_{b\in\b}\x_bb$.
Let $V_{\ge0}={}^\l V_{\ge0}$ (resp. $V_\RR$) be the set of all $\x\in V$ such that $\x_b\in\RR_{\ge0}$ (resp. 
$\x_b\in\RR$) for any $b\in\b$.  We have $V_{\ge0}\sub V_\RR$. Note that $V_\RR$ is stable under the action 
of $G(\RR)$ on $V$. Let $P_{\ge0}={}^\l P_{\ge0}$ (resp. $P_\RR$) be the set of lines $L\in P$ such that 
$L\cap V_{\ge0}\ne0$ (resp. $L\cap V_\RR\ne0$.) We have $P_{\ge0}\sub P_\RR$.

Let $V^\bul_{\ge0}={}^\l V^\bul_{\ge0}=V^\bul\cap V_{\ge0}$, 
$P^\bul_{\ge0}={}^\l P^\bul_{\ge0}=P^\bul\cap P_{\ge0}$. 

\mpb

Now let $K$ be a semifield. Let $V(K)={}^\l V(K)$ be the set of formal sums 
$\x=\sum_{b\in\b}\x_bb,\x_b\in K^!$. This 
is a monoid under addition $(\sum_{b\in\b}\x_bb)+(\sum_{b\in\b}\x'_bb)=\sum_{b\in\b}(\x_b+\x'_b)b$ and we 
define scalar multiplication $K^!\T V(K)@>>>V(K)$ by $(k,\sum_{b\in\b}\x_bb)\m\sum_{b\in\b}(k\x_b)b$. 

For $\x=\sum_{b\in\b}\x_bb\in V(K)$ we define $\supp(\x)=\{b\in\b;\x_b\in K\}$. 

Let $\End(V(K))$ be the set of maps $\z:V(K)@>>>V(K)$ such that $\z(\x+\x')=\z(\x)+\z(\x')$ for 
$\x,\x'$ in $V(K)$ and $\z(k\x)=k\z(\x)$ for $\x\in V(K),k\in K^!$. This is a monoid under composition
of maps. Define $\un\circ\in V(K)$ by $\un\circ_b=\circ$ for all $b\in\b$. The group $K$ 
(for multiplication in the semifield structure) acts freely (by scalar multiplication) on $V(K)-\un\circ$; 
let $P(K)={}^\l P(K)$ be the set of orbits of this action.

For $i\in I,n\in\ZZ$ we define $e_i^{(n)},f_i^{(n)}$ in $\End(V(K))$ by
$$e_i^{(n)}(b)=\sum_{b'\in\b}c_{b,b',i,n}b',\qua 
f_i^{(n)}(b)=\sum_{b'\in\b}d_{b,b',i,n}b',$$
with $b\in\b$. Here a natural number $N$ (such as $c_{b,b',i,n}$ or $d_{b,b',i,n}$)
is viewed as an element of $K^!$ given by $1+1+\do+1$ ($N$ terms, where $1$ is the neutral element for the
product in $K$, if $N>0$) or by $\circ\in K^!$ (if $N=0$).

For $i\in I,k\in K$ we define $i^k\in\End(V(K))$, $(-i)^k\in\End(V(K))$ by
$$i^k(b)=\sum_{n\in\NN}k^ne_i^{(n)}b,\qua (-i)^k(b)=\sum_{n\in\NN}k^nf_i^{(n)}b,$$ 
for any $b\in\b$. We show:

(a) {\it The map $i^k:V(K)@>>>V(K)$ is injective. The map $(-i)^k:V(K)@>>>V(K)$ is injective.}
\nl
Using a partial order of the weights of $V$, we can write $V(K)$ as a direct sum of monoids
$V(K)_s,s\in\ZZ$ where $V(K)_s=\{\un\circ\}$ for all but finitely many $s$ and $(-i)^k$ maps any
$\x\in V(K)_s$ to $\x$ plus an element in the direct sum of $V(K)_{s'}$ with $s'<s$.
Then (a) for $(-i)^k$ follows immediatly. A similar proof applies to $i^k$.

\mpb

For $i\in I,k\in K$ we define $\un i^k\in\End(V(K))$ by $\un i^k(b)=k^{z_i(b)}b$ for any $b\in\b$.
Let $\fG(K)$ be the monoid associated to $G,K$ by generators and relations in \cite{L19a, 2.10(i)-(vii)}.
(In {\it loc.cit.} it is assumed that $K$ is as in 0.1(i) or 0.1(ii) but the same definition makes sense
for any $K$.) We have the following result.

\proclaim{Proposition 1.5} The elements $i^k,(-i)^k,\un i^k$ (with $i\in I,k\in K$) in $\End(V(K))$
satisfy the relations in \cite{L19a, 2.10(i)-(vii)} defining the monoid $\fG(K)$ hence they define a monoid 
homomorphism $\fG(K)@>>>\End(V(K))$.
\endproclaim
We write the relations in {\it loc.cit.} (for the semifield $\RR_{>0}$)
for the endomorphisms $x_i(t),y_i(t),i(t)$ of $V$ with $t\in R_{>0}$. These relations can be expressed as
a set of identities satisfied by $c_{b,b',i,n}$, $d_{b,b',i,n}$, $z_i(b)$ and these identities show that
the endomorphisms $i^k,(-i)^k,\un i^k$ of $V(K)$ satisfy the relations in {\it loc.cit.} (for the 
semifield $K$). The result follows.

\subhead 1.6\endsubhead
Consider a homomorphism of semifields $r:K_1@>>>K_2$. Now $r$ induces a homomorphism of monoids
$\fG_r:\fG(K_1)@>>>\fG(K_2)$. It also induces a homomorphism of monoids $V_r:V(K_1)@>>>V(K_2)$ given by
$\sum_{b\in\b}\x_bb\m\sum_{b\in\b}r(\x_b)b$. From the definitions, for $g\in\fG(K_1),\x\in V(K_1)$, we 
have $V_r(g\x)=\fG_r(g)(V_r(\x))$ where $g\x$ is given by the $\fG(K_1)$-action on $V(K_1)$ and 
$\fG_r(g)(V_r(\x))$ is given by the $\fG(K_2)$-action on $V(K_2)$.
Assuming that $r:K_1@>>>K_2$ is surjective (so that $\fG_r:\fG(K_1)@>>>\fG(K_2)$ is surjective) we deduce:

(a) {\it If $E$ is a subset of $V(K_1)$ which is stable under the $\fG(K_1)$-action on $V(K_1)$, then
the subset $V_r(E)$ of $V(K_2)$ is stable under the $\fG(K_2)$-action on $V(K_2)$.}

\subhead 1.7\endsubhead
In the remainder of this section we assume that $\l\in\cx^{++}$. Then $L\m G_L$ is an isomorphism 
$\p:P^\bul@>\si>>\cb$ and  

(a) $\p$ restricts to a bijection $\p_{\ge0}:P^\bul_{\ge0}@>\si>>\cb_{\ge0}$.
\nl
See \cite{L94b, 8.17}.

\subhead 1.8\endsubhead
Let $\Om$ be the set of all open nonempty subsets of $\CC$. Let $X$ be an algebraic variety over $\CC$.
Let $X_1$ be the set of pairs $(U,f_U)$ where $U\in\Om$ and $f_U:U@>>>X$ is a morphism of algebraic varieties.
We define an equivalence relation on $X_1$ in which $(U,f_U),(U',f_{U'})$ are equivalent if
$f_U|_{U\cap U'}=f_{U'}|_{U\cap U'}$. Let $\tX$ be the set of equivalence classes. An element of $\tX$ is 
said to be a rational map $f:\CC\rhd X$. For $f\in\tX$ let $\Om_f$ be the set of all $U\in\Om$ such that $f$ 
contains $(U,f_U)\in X_1$ for some $f_U$; we shall then write $f(t)=f_U(t)$ for $t\in U$.
We shall identify any $x\in X$ with the constant map $f_x:\CC@>>>X$ with image $\{x\}$; thus $X$ can be 
identified with a subset of $\tX$.
If $X'$ is another algebraic variety over $\CC$ then we have $\wt{X\T X'}=\tX\T\tX'$ canonically.
If $F:X@>>>X'$ is a morphism then there is an induced map $\tF:\tX@>>>\tX'$; to $f:\CC\rhd X$ it attaches 
$f':\CC\rhd X'$ where for some $U\in\Om_f$ we have $f'(t)=F(f(t))$ for all $t\in U$. 
If $H$ is an algebraic group over $\CC$ then $\tH$ is a group with multiplication 
$\tH\T\tH=\wt{H\T H}@>>>\tH$ induced by the multiplication map $H\T H@>>>H$. Note that 
$H$ is a subgroup of $\tH$. In particular, the group $\tG$ is defined. Also, the additive group $\ti\CC$ 
and the multiplicative group $\wt{\CC^*}$ are defined. Also $\tcb$ is defined.

\subhead 1.9\endsubhead
Let $X$ be an algebraic variety over $\CC$ with a given subset $X_{\ge0}$. We define a subset $\tX_{\ge0}$ 
of $\tX$ as follows: $\tX_{\ge0}$ is the set of all $f\in\tX$ such that for some $U\in\Om_f$ and some 
$\e\in\RR_{>0}$ we have $(0,\e)\sub U$ and $f(t)\in X_{\ge0}$ for all $t\in(0,\e)$.
(In particular, $\tG_{\ge0}$ is defined in terms of $G,G_{\ge0}$
and $\tcb_{\ge0}$ is defined in terms of $\cb,\cb_{\ge0}$.)
If $X'$ is another algebraic variety over $\CC$ with a given subset $X'_{\ge0}$, then $X\T X'$ with its
subset $(X\T X')_{\ge0}=X_{\ge0}\T X'_{\ge0}$ gives rise as above to the set $\wt{X\T X'}_{\ge0}$ which can 
be identified with $\tX_{\ge0}\T\tX'_{\ge0}$.
If $F:X@>>>X'$ is a morphism such that $F(X_{\ge0})\sub X'_{\ge0}$, then the induced map
$\tF:\tX@>>>\tX'$ carries $\tX_{\ge0}$ into $\tX'_{\ge0}$ hence it restricts to a map
$\tF_{\ge0}:\tX_{\ge0}@>>>\tX'_{\ge0}$. From the definitions we see that:

(a) {\it if $\tF$ is an isomorphism of $\tX$ onto an open subset of $\tX'$ and $F$ carries
$\tX_{\ge0}$ bijectively onto $\tX'_{\ge0}$, then the map $\tF_{\ge0}$ is a bijection.}
\nl
Now the multiplication $G\T G@>>>G$ carries $G_{\ge0}\T G_{\ge0}$ to $G_{\ge0}$ hence it 
induces a map $\tG_{\ge0}\T\tG_{\ge0}@>>>\tG_{\ge0}$ which makes $\tG_{\ge0}$ into a monoid; the conjugation 
action $G\T\cb@>>>\cb$ carries $G_{\ge0}\T\cb_{\ge0}$ to $\cb_{\ge0}$ hence it induces a map 
$\tG_{\ge0}\T\tcb_{\ge0}@>>>\tcb_{\ge0}$ which define an action of the monoid $\tG_{\ge0}$ on $\tcb_{\ge0}$.
We define $\ti{\CC^*}_{\ge0}$ in terms of $\CC^*$ and its subset $\CC^*_{\ge0}:=\RR_{>0}$. The multiplication 
on $\CC^*$ preserves $\CC^*_{\ge0}$ hence it induces a map 
$\ti{\CC^*}_{\ge0}\T\ti{\CC^*}_{\ge0}@>>>\ti{\CC^*}_{\ge0}$ which makes $\ti{\CC^*}_{\ge0}$ into an abelian 
group. 
We define $\ti{\CC}_{\ge0}$ in terms of $\CC$ and its subset $\CC_{\ge0}:=\RR_{\ge0}$. The addition 
on $\CC$ preserves $\CC_{\ge0}$ hence it induces a map 
$\ti{\CC}_{\ge0}\T\ti{\CC}_{\ge0}@>>>\ti{\CC}_{\ge0}$ which makes $\ti{\CC}_{\ge0}$ into an abelian monoid.
The imbedding $\CC^*\sub\CC$ induces an imbedding $\ti{\CC^*}_{\ge0}@>>>\ti{\CC}_{\ge0}$; the monoid 
operation on
$\ti{\CC}_{\ge0}$ preserves the subset $\ti{\CC^*}_{\ge0}$ and makes $\ti{\CC^*}_{\ge0}$ into an abelian
monoid. This, together with the multiplication on $\ti{\CC^*}_{\ge0}$ makes
$\ti{\CC^*}_{\ge0}$ into a semifield. From the definitions we see that this semifield is the same as
$K$ in 0.1(i) and that $\tG_{\ge0}$ is the monoid  associated to $G$ and $K$ in \cite{L94b, 2.2} (which is the
sme as $\fG(K)$). We define $\cb(K)$ to be $\tcb_{\ge0}$ with the action of $\tG_{\ge0}=\fG(K)$ described above.
This achieves what was stated in 0.1 for $K$ as in 0.1(i).

\subhead 1.10\endsubhead
In the remainder of this section $K$ will denote the semifield in 0.1(i) and we assume that $\l\in\cx^{++}$. 
We associate $\tP_{\ge0}={}^\l\tP_{\ge0}$ to $P$ and its subset $P_{\ge0}$ as in 1.9. 
We associate $\tP^\bul_{\ge0}={}^\l\tP^\bul_{\ge0}$ to 
$P^\bul$ and its subset $P^\bul_{\ge0}$ as in 1.9. We write
$P^\bul(K)={}^\l P^\bul(K)=\tP^\bul_{\ge0}$. 

\mpb

We associate $\tV_{\ge0}={}^\l\tV_{\ge0}$ 
to $V$ and its subset $V_{\ge0}$ as in 1.9. We can identify $\tV_{\ge0}=V(K)$ 
(see 1.4). We associate $\tV^\bul_{\ge0}={}^\l\tV^\bul_{\ge0}$ to 
$V^\bul$ and its subset $V^\bul_{\ge0}$ as in 1.9. 
We write $V^\bul(K)={}^\l V^\bul(K)=\tV^\bul_{\ge0}$. We have $V^\bul(K)\sub\tV_{\ge0}$.

The obvious map $a':V-0@>>>P$ restricts to a (surjective) map 
$a'_{\ge0}:V_{\ge0}-0@>>>P_{\ge0}$ and defines a map $\ta'_{\ge0}:\tV_{\ge0}-0@>>>\tP_{\ge0}$. 
The scalar multiplication $\CC^*\T(V-0)@>>>V-0$ carries $\CC^*_{\ge0}\T(V_{\ge0}-0)$ to $V_{\ge0}-0$ hence it
induces a map $\wt{\CC^*}_{\ge0}\T(\tV_{\ge0}-0)@>>>\tV_{\ge0}-0$ which is 
a (free) action of the group $K=\wt{\CC^*}_{\ge0}$ on $\tV_{\ge0}-0=V(K)-0$. 
From the definitions we see that $\ta'_{\ge0}$ is surjective and it induces a bijection 
$(V(K)-0)/K@>\si>>\tP_{\ge0}$. Thus we have $\tP_{\ge0}=P(K)$ (notation of 1.4). Note that $P^\bul(K)\sub P(K)$.

The obvious map $a:V^\bul-0@>>>P^\bul$ restricts to a (surjective) map 
$a_{\ge0}:V^\bul_{\ge0}-0@>>>P^\bul_{\ge0}$ 
and it defines a map $\ta_{\ge0}:V^\bul(K)=\tV^\bul_{\ge0}-0@>>>\tP^\bul_{\ge0}=P^\bul(K)$. 
The (free) $K$-action on $\tV_{\ge0}-0$ considered above restricts to a (free) $K$-action
on $V^\bul(K)-0=\tV^\bul_{\ge0}-0$. From the definitions we see that
$\ta_{\ge0}$ is constant on any orbit of this action. We show:

(a) {\it The map $\ta_{\ge0}$ is surjective. It induces a bijection $(V^\bul(K)-0)/K@>\si>>P^\bul(K)$.}
\nl
Let $f\in\tP^\bul_{\ge0}$. We can find $U\in\Om_f$, $\e\in\RR_{>0}$ such that
$(0,\e)\sub U$ and $f(t)\in P^\bul_{\ge0}$ for $t\in(0,\e)$. Using the surjectivity of $a_{\ge0}$ we see that
for $t\in(0,\e)$ we have $f(t)=a(x_t)$ where $t\m x_t$ is a function $(0,\e)@>>>V^\bul_{\ge0}-0$.
We can assume that there exists $B\in\cb(\RR)$ such that 
$\p(f(t))$ is opposed to $B$ for all $t\in U$. Let
$\co=\{B_1\in\cb;B_1\text{ opposed to }B\}$; thus we have $\p(f(t))\in\co$ for all $t\in U$. Let 
$B'\in\co\cap\cb(\RR)$ and let $\x'\in V_\RR-0$ be such that $\p(\CC\x')=B'$.
Let $U_B$ be the unipotent radical of $B$. Then $U_B@>>>\co$, $u\m uB'u\i$ is an isomorphism. Hence 
there is a unique morphism $\z:\co@>>>V^\bul-0$ such 
that $\z(uB'u\i)=u\x'$ for any $u\in U_B$. From the definitions we have 
$\z(\co\cap\cb(\RR))\sub(V_\RR\cap V^\bul)-0$. We define $f':U@>>>V^\bul-0$ by $f'(t)=\z(\p(f(t)))$.
We can view $f'$ as an element of $\tV^\bul-0$ such that $\ta(f')=f$. Since $\p(f(t))\in\cb(\RR)$, we have
$f'(t)\in (V_\RR\cap V^\bul)-0$ for $t\in(0,\e)$. For such $t$ we have $a(f'(t))=f(t)=a(x_t)$ hence
$f'(t)=z_tx_t$ where $t\m z_t$ is a (possibly discontinuous) function $(0,\e)@>>>\RR-0$. 
Since $x_t\in V_{\ge0}-0$ and $\RR_{>0}(V_{\ge0}-0)=V_{\ge0}-0$, we see that for $t\in(0,\e)$ we have
$f'(t)\in(V_{\ge0}-0)\cup(-1)(V_{\ge0}-0)$. Since $(0,\e)$ is connected and $f'$ is continuous
(in the standard topology) we see that $f'(0,\e)$ is contained in one of the connected components of
$(V_{\ge0}-0)\cup(-1)(V_{\ge0}-0)$ that is, in either $V_{\ge0}-0$ or in $(-1)(V_{\ge0}-0)$. Thus there
exists $s\in\{1,-1\}$ such that $sf'(0,\e)\sub V_{\ge0}-0$ hence also $sf'(0,\e)\sub V^\bul_{\ge0}-0$.
We define $f'':U@>>>V^\bul-0$ by $f''(t)=sf'(t)$. We can view $f''$ as an element of $\tV^\bul_{\ge0}-0$ such 
that $\ta_{\ge0}(f')=f$. This proves that $\ta_{\ge0}$ is surjective.
The remaining statement of (a) is immediate.

\mpb

Since $P^\bul$ and its subset $P^\bul_{\ge0}$ can be identified with $\cb$ and its 
subset $\cb_{\ge0}$ (see 1.7(a)), we see that we may identify $P^\bul(K)=\cb(K)$.
The action of $\fG(K)$ on $P^\bul(K)$ induced from that on $V^\bul(K)-0$ is the same as the previous action of 
$\fG(K)$, see \cite{L19a, 2.13(d)}. This gives a second incarnation of $\cb(K)$.

\subhead 1.11\endsubhead
Let $\ZZ$ be the semifield in 0.1(ii). Following \cite{L94b}, we define 
a (surjective) semifield homomorphism $r:K@>>>\ZZ$ by $r(x^ef_1/f_2)=e$ (notation of 0.1). Now $r$ induces
a surjective map $V_r:V(K)@>>>V(\ZZ)$ as in 1.6. Let $V^\bul(\ZZ)={}^\l V^\bul(\ZZ)\sub V(\ZZ)$ 
be the image under $V_r$ of the subset $V^\bul(K)$ of $V(K)$. Then $V^\bul(\ZZ)-\un\circ=V_r(V^\bul(K)-0)$.

The $\ZZ$-action on $V(\ZZ)-\un\circ$ in 1.4 leaves $V^\bul(\ZZ)-\un\circ$ stable. (We use the $K$-action on 
$V^\bul(K)-0$.) 
Let $P^\bul(\ZZ)={}^\l P^\bul(\ZZ)$ 
be the set of orbits of this action.  We have $P^\bul(\ZZ)\sub P(\ZZ)$ (notation of 1.4).
From 1.6(a) we see that $V^\bul(\ZZ)-\un\circ$ is stable under the $\fG(\ZZ)$-action on $V(\ZZ)$ in
1.6. Since the $\fG(\ZZ)$-action commutes with scalar multiplication by $\ZZ$ it follows that
the $\fG(\ZZ)$-action on $V(\ZZ)-\un\circ$ and $V^\bul(\ZZ)-\un\circ$ induces a $\fG(\ZZ)$-action on 
$P(\ZZ)$ and $P^\bul(\ZZ)$. 

\subhead 1.12\endsubhead
We set $\cb(\ZZ)={}^\l P^\bul(\ZZ)$. This achieves what was stated in 0.1 for the semifield $\ZZ$. 
This definition of $\cb(\ZZ)$ depends on the choice of 
$\l\in\cx^{++}$. In \S4 we will show that $\cb(\ZZ)$ is independent of this choice up to a canonical bijection.
(Alternatively, if one wants a definition without such a choice one could take $\l$ such that
$\la i,\l\ra=1$ for all $i\in I$.)

\head 2. Preparatory results\endhead
\subhead 2.1\endsubhead
We preserve the setup of 1.4. As shown in \cite{L94a, 5.3, 4.2}, for $w\in W$ and  
$\ii=(i_1,i_2,\do,i_m)\in\ci_w$, the subspace of $V$ generated by the vectors 
$$f_{i_1}^{(c_1)}f_{i_2}^{(c_2)}\do f_{i_m}^{(c_m)}\x^+$$ 
for various $c_1,c_2,\do,c_m$ in $\NN$ is independent of $\ii$ (we denote it by $V^w$) and $\b^w:=\b\cap V^w$ is 
a basis of it. Let $V'{}^\ii$ be the subspace of $V$ generated by the vectors
$$e_{i_m}^{(d_m)}e_{i_{m-1}}^{(d_{m-1})}\do e_{i_1}^{(d_1)}b_w$$
for various $d_1,d_2,\do,d_m$ in $\NN$, where 
$$b_w=\dw\x^+,$$ 
$$\dw=\ds_{i_1}\ds_{i_2}\do\ds_{i_m}.$$
We show:
$$V^w=V'{}^\ii.\tag a$$
We show that $V^w\sub V'{}^\ii$. We argue by induction on $m=|w|$. If $m=0$, the result is obvious. Assume now 
that $m\ge1$. Let $c_1,c_2,\do,c_m$ be in $\NN$. By the induction hypothesis, 
$$f_{i_1}^{(c_1)}f_{i_2}^{(c_2)}\do f_{i_m}^{(c_m)}\x^+\tag b$$
is a linear combination of vectors of form
$$f_{i_1}^{(c_1)}e_{i_m}^{(d_m)}e_{i_{m-1}}^{(d_{m-1})}\do e_{i_2}^{(d_2)}b_{s_{i_1}w}$$
for various $d_2,\do,d_m$ in $\NN$. 
Using the known commutation relations between $f_{i_1}$ and $e_j$ we see that (b) is a linear combination of 
vectors of form
$$e_{i_m}^{(d_m)}e_{i_{m-1}}^{(d_{m-1})}\do e_{i_2}^{(d_2)}f_{i_1}^{(c_1)}b_{s_{i_1}w}$$
for various $d_2,\do,d_m$ in $\NN$. It is then enough to show that
$$f_{i_1}^{(c_1)}b_{s_{i_1}w}=e_{i_1}^{(d_1)}\ds_{i_1}b_{s_{i_1}w}$$
for some $d_1\in\NN$. This follows from the fact that 

(c) {\it $e_{i_1}b_{s_{i_1}w}=0$ and $b_{s_{i_1}w}$ is in a weight space of $V$.}
\nl
Next we show that $V'{}^\ii\sub V^w$. We argue by induction on $m=|w|$. If $m=0$ the result is obvious. Assume 
now that $m\ge1$. Since $V^w$ is stable under the action of $e_i (i\in I)$, it is enough to show that $b_w\in V^w$. 
By the induction hypothesis, $b_{s_{i_1}w}\in V^{s_{i_1w}}$. Using (c), we see that for some $c_1\in\NN$ we have 
$$b_w=\ds_{i_1}b_{s_{i_1}w}=f_{i_1}^{(c_1)}b_{s_{i_1}w}\in f_{i_1}^{(c_1)}V^{s_{i_1}w}\sub V^w.$$
This completes the proof of (a).

\mpb

From \cite{L93, 28.1.4} one can deduce that $b_w\in\b$. From (a) we see that $b_w\in V^w$. It follows that
$$b_w\in\b^w.\tag d$$

\subhead 2.2\endsubhead
For $v\le w$ in $W$ we set 
$$\cb_{v,w}=\{B\in\cb,pos(B^+,B)=w,pos(B^-,B)=w_Iv\}$$ 
(a locally closed subvariety of $\cb$) and
$$(\cb_{v,w})_{\ge0}=\cb_{\ge0}\cap\cb_{v,w}.$$ 
We have $\cb=\sqc_{v\le w\text{ in }W}\cb_{v,w}$, $\cb_{\ge0}=\sqc_{v\le w\text{ in }W}(\cb_{v,w})_{\ge0}$.

\subhead 2.3\endsubhead
Recall that there is a unique isomorphism $\ph:G@>>>G$ such that $\ph(x_i(t))=y_i(t)$, $\ph(y_i(t))=x_i(t)$ 
for all $i\in I,t\in\CC$ and $\ph(g)=g\i$ for all $g\in T$. This carries Borel subgroups to Borel subgroups
hence induces an isomorphism $\ph:\cb@>>>\cb$ such that $\ph(B^+)=B^-$, $\ph(B^-)=B^+$. For $i\in I$ we have 
$\ph(\ds_i)=\ds_i\i$. Hence $\ph$ induces the identity map on $W$. For $v\le w$ in $W$ we have 
$ww_I\le vw_I$; moreover,

(a) {\it $\ph$ defines an isomorphism $\cb_{ww_I,vw_I}@>\si>>\cb_{v,w}$.}
\nl
(See \cite{L19b, 1.4(a)}.) From the definition we have 

(b) $\ph(G_{\ge0})=G_{\ge0}$.
\nl
From \cite{L94b, 8.7} it follows that 

(c) $\ph(\cb_{\ge0})=\cb_{\ge0}$.
\nl
From (a),(c) we deduce:

(d) {\it $\ph$ defines a bijection $(\cb_{ww_I,vw_I})_{\ge0}@>\si>>(\cb_{v,w})_{\ge0}$.}
\nl
By \cite{L90b, \S3} there is a unique linear isomorphism $\ph:V@>>>V$ such that
$\ph(g\x)=\ph(g)\ph(\x)$ for all $g\in G$, $\x\in V$ and such that $\ph(\x^+)=\x^-$; we have
$\ph(\b)=\b$ and $\ph^2(\x)=\x$ for all $\x\in V$.

\subhead 2.4\endsubhead
Assume now that $\l\in\cx^{++}$.
Let $B\in\cb_{v,w}$ and let $L\in P^\bul$ be such that $\p(L)=B$. Let $\x\in L-0$, $b\in\b$. We show:

(a) $\x_b\ne0\implies b\in\b^w\cap\ph(\b^{vw_I})$.
\nl
We have $B=gB^+g\i$ for some $g\in B^+\dw B^+$. Then $\x=cg\x^+$ for some $c\in\CC^*$. We write 
$g=g'\dw g''$ with $g'\in U^+,g''\in B^+$. We have $\x=c'g'\dw\x^+=c'g'b_w$ where $c'\in\CC^*$.
By 2.1(d) we have $b_w\in\b^w$. Moreover, $V^w$ is stable by the action of $U^+$; we see that $\x\in V^w$.
Since $\x_b\ne0$ we have $b\in\b^w$.
Let $B'=\ph(B)$. We have $B'\in\cb_{ww_I,vw_I}$ (see 2.3(a)). Let $L'=\ph(L)\in P^\bul$ and let
$\x'=\ph(\x)\in L'-0$, $b'=\ph(b)\in\b$. We have $\x'_{b'}\ne0$. Applying the first part of the proof with 
$B,L,\x,v,w,b$ replaced by $B',L',\x',v',w',b'$ we obtain $b'\in\b^{vw_I}$. Hence $b\in\ph(\b^{vw_I})$. Thus, 
$b\in\b^w\cap\ph(\b^{vw_I})$, as required.

\subhead 2.5\endsubhead
We return to the setup of 1.4. For $i\in I$ we set 
$$V^{e_i}=\{\x\in V;e_i(\x)=0\}=\{\x\in V;\sum_{b\in\b}\x_bc_{b,b',i,1}=0 \text{ for all }b'\in\b\},$$
$$V^{f_i}=\{\x\in V;f_i(\x)=0\}=\{\x\in V;\sum_{b\in\b}\x_bd_{b,b',i,1}=0 \text{ for all }b'\in\b\}.$$
If $\x\in V_{\ge0}$, the condition that $\sum_{b\in\b}\x_bc_{b,b',i,1}=0$ is equivalent to the condition 
that $\x_bc_{b,b',i,1}=0$ for any $b,b'$ in $\b$. Thus we have
$$V_{\ge0}\cap V^{e_i}=\{\x\in V_{\ge0}; \x=\sum_{b\in\b^{e_i}}\x_bb\}$$
where $\b^{e_i}=\{b\in\b;c_{b,b',i,1}=0\text{ for any }b'\in\b\}$. Similarly, we have
$$V_{\ge0}\cap V^{f_i}=\{\x\in V_{\ge0}; \x=\sum_{b\in\b^{f_i}}\x_bb\}$$
where $\b^{f_i}=\{b\in\b;d_{b,b',i,1}=0\text{ for any }b'\in\b\}$.

Now the action of $\ds_i$ on $V$ defines an isomorphism $\ct_i:V^{e_i}@>>>V^{f_i}$.
If $b\in\b^{e_i}$ we have $\ct_i(b)=f_i^{(\la i,\nu_b\ra)}b=\sum_{b'\in\b}d_{b,b',i,\la i,\nu_b\ra}b'$; 
in particular, we 
have $\ct_i(b)\in V_{\ge0}\cap V^{f_i}$. Thus $\ct_i$ restricts to a map 
$\ct'_i:V_{\ge0}\cap V^{e_i}@>>>V_{\ge0}\cap V^{f_i}$. Similarly the action of $\ds_i\i$ restricts to a map 
$\ct''_i:V_{\ge0}\cap V^{f_i}@>>>V_{\ge0}\cap V^{e_i}$. This is clearly the inverse of $\ct'_i$.

\subhead 2.6\endsubhead
Now let $K$ be a semifield. Let 
$$V(K)^{e_i}=\{\sum_{b\in\b}\x_bb;\x_b\in K^!\text{ if }b\in\b^{e_i},\x_b=\circ\text{ if }b\in\b-\b^{e_i}\},$$
$$V(K)^{f_i}=\{\sum_{b\in\b}\x_bb;\x_b\in K^!\text{ if }b\in\b^{f_i},\x_b=\circ\text{ if }b\in\b-\b^{f_i}\}.$$
We define $\ct_{i,K}:V(K)@>>>V(K)$ by
$$\sum_{b\in\b}\x_bb\m\sum_{b'\in\b}(\sum_{b\in\b}d_{b,b',i,\la i,\nu_b\ra}\x_b)b'$$
(notation of 1.4). From the results in 2.5 one can deduce that 

(a) {\it $\ct_{i,K}$ restricts to a bijection $\ct'_{i,K}:V(K)^{e_i}@>\si>>V(K)^{f_i}$.}

\subhead 2.7\endsubhead
Let $K$ be a semifield. We define an involution  $\ph:V(K)@>>>V(K)$ by 
$\ph(\sum_{b\in\b}\x_bb)=\sum_{b\in\b}\x_{\ph(b)}b$. (Here $\x_b\in K^!$; we use that $\ph(\b)=\b$.) This 
restricts to an involution $V(K)-\un\circ@>>>V(K)-\un\circ$ which induces an involution $P(K)@>>>P(K)$ denoted 
again by $\ph$.

\head 3. Parametrizations\endhead
\subhead 3.1\endsubhead
In this section $K$ denotes the semifield in 0.1(i). For $v\le w$ in $W$ we define 
$\cb_{v,w}(K)=\wt{\cb_{v,w}}_{\ge0}$ as in 1.9 in terms of $\cb_{v,w}$ and its subset 
$(\cb_{v,w})_{\ge0}$. We have 
$$\cb(K)=\sqc_{v\le w\text{ in }W}\cb_{v,w}(K).$$

\subhead 3.2\endsubhead
We preserve the setup of 1.4. 
We now fix $v\le w$ in $W$ and $\ii=(i_1,i_2,\do,i_m)\in\ci_w$. According to \cite{MR},
there is a unique sequence $q_1,q_2,\do,q_m$  with $q_k\in\{s_{i_k},1\}$
for $k\in[1,m]$, $q_1q_2\do q_m=v$ and such that 
$q_1\le q_1q_2\le\do\le q_1q_2\do q_m$ and $q_1\le q_1s_{i_2},q_1q_2\le q_1q_2s_{i_3},\do,
q_1q_2\do q_{m-1}\le q_1q_2\do q_{m-1}s_{i_m}$. Let $[1,m]'=\{k\in[1,m];q_k=1\}$, 
$[1,m]''=\{k\in[1,m];q_k=s_{i_k}\}$. 
Let $A$ be the set of maps $h:[1,m]'@>>>\CC^*$; this is naturally an algebraic variety over $\CC$.
Let $A_{\ge0}$ be the subset of $A$ consisting of maps $h:[1,m]'@>>>\RR_{>0}$. 
Following \cite{MR}, we define a morphism $\s:A@>>>G$ by $h\m g(h)_1g(h)_2\do g(h)_m$ where 

(a) {\it $g(h)_k=y_{i_k}(h(k))$ if $k\in[1,m]'$ and $g(h)_k=\ds_{i_k}$ if $k\in[1,m]''$.}
\nl
We show:

(b) {\it If $h\in A_{\ge0}$, then $\s(h)\x^+\in V^w$, so that $\s(h)$ is a 
linear combination of vectors $b\in\b^w$. Moreover, $(\s(h)\x^+)_{b_w}\ne0$.}
\nl
From the properties of Bruhat decomposition, for any $h\in A_{\ge0}$ we have $\s(h)\in B^+\dw B^+$, so that 
$\s(h)\x^+=cu\dw\x^+=cub_w$ where $c\in\CC^+$, $u\in U^+$. Since $b_w\in V^w$ and $V^w$ is stable under the 
action of $U^+$, it follows that $cu\dw\x^+\in V^w$. More precisely, $ub_w=b_w$ plus a linear combination of 
elements $b\in\b$ of weight other than that of $b_w$. This proves (b).

We show:

(c) {\it Let $h\in A_{\ge0}$. Assume that $i\in I$ is such that $|s_iw|>|w|$ and that $b\in\b$ is such that
$(\s(h)\x^+)_b\ne0$. Then $\nu_b\ne\nu_{b_w}+i'$.}
\nl
Since $|s_iw|>|w|$ we have $e_ib_w=0$. We write $\s(h)x^+=cub_w$ with $c,u$ as in the proof of (b). Now $ub_w$ 
is a linear combination of vectors of the form $e_{j_1}e_{j_2}\do e_{j_k}b_w$ with $j_t\in I$.
Such a vector is in a weight space $V(\nu)$ with $\nu=\nu_{b_w}+j'_1+j'_2+\do+j'_k$. If 
$j'_1+j'_2+\do+j'_k=i'$ then $k=1$ and $j_1=i$. But in this case we have $e_{j_1}e_{j_2}\do e_{j_k}b_w=e_ib_w=0$.
The result follows. 

\subhead 3.3\endsubhead 
Let $h\in A_{\ge0}$. Let $k\in[1,m]''$. The following result appears in the proof of \cite{MR, 11.9}.

(a) {\it We have $(g(h)_{k+1}g(h)_{k+2}\do g(h)_m)\i x_{i_k}(a)g(h)_{k+1}g(h)_{k+2}\do g(h)_m\in U^+$.}
\nl
From (a) it follows that for $\x\in V$ we have 
$$e_{i_k}(g(h)_{k+1}g(h)_{k+2}\do g(h)_m\x)=g(h)_{k+1}g(h)_{k+2}\do g(h)_m (e'\x)$$
where $e':V@>>>V$ is a linear combination of products of one or more factors $e_j,j\in I$.
When $\x=\x^+$ we have $e'\x=0$ hence $e_{i_k}(g(h)_{k+1}g(h)_{k+2}\do g(h)_m\x^+)=0$. We can write uniquely 
$$g(h)_{k+1}g(h)_{k+2}\do g(h)_m\x^+=\sum_{\nu\in\cx}(g(h)_{k+1}g(h)_{k+2}\do g(h)_m\x^+)_\nu$$
with $(g(h)_{k+1}g(h)_{k+2}\do g(h)_m\x^+)_\nu\in V_\nu$. We have 
$$\sum_{\nu\in\cx}e_{i_k}((g(h)_{k+1}g(h)_{k+2}\do g(h)_m \x^+)_\nu)=0.$$
Since the elements $e_{i_k}((g(h)_{k+1}g(h)_{k+2}\do g(h)_m \x^+)_\nu)$ (for various $\nu\in\cx$) are 
in distinct weight spaces, it follows that $e_{i_k}((g(h)_{k+1}g(h)_{k+2}\do g(h)_m \x^+)_\nu)=0$ for any 
$\nu\in\cx$. If $\x\in V_\nu$ satisfies $e_{i_k}\x=0$, then 

(b) $\ds_{i_k}\x=f_{i_k}^{(\la i_k,\nu\ra)}\x$.
\nl
(If $\la i_k,\nu\ra<0$ then $\x=0$ so that both sides of (b) are $0$.) We deduce
$$g(h)_k((g(h)_{k+1}g(h)_{k+2}\do g(h)_m \x^+)_\nu)=
f_{i_k}^{(\la i_k,\nu\ra)}((g(h)_{k+1}g(h)_{k+2}\do g(h)_m \x^+)_\nu)\tag c$$
for any $\nu\in\cx$.

\subhead 3.4\endsubhead 
Let $h\in A_{\ge0}$. For any $k\in[1,m]$ we set $[k,m]'=[k,m]\cap[1,m]'$, $[k,m]''=[k,m]\cap[1,m]''$.
Let $\ce_{\ge k}$ be the set of all maps $\c:[k,m]'@>>>\NN$. (If $[k,m]'=\emp$, $\ce_{\ge k}$ consists of a 
single element.) 
For $\c\in\ce_{\ge k}$ and $k'\in[k,m]$ let $\c_{\ge k'}$ be the restriction of $\c$ to $[k',m]'$.

We now define an integer $c(k,\c)$ for any $k\in[1,m]''$ and any $\c\in\ce_{\ge k}$ by descending induction 
on $k$. We can assume that $c(k',\c')$ is defined for any $k'\in[k+1,m]''$ and any $\c'\in\ce_{\ge k'}$. We set 
$c_{k,\c}=\la i_k,\nu\ra$ where
$$\nu=\l-\sum_{\k\in[k+1,m]'}\c(\k)i'_\k-\sum_{\k\in[k+1,m]'';c(\k,\c_{\ge\k})\ge0}c(\k,\c_{\ge\k})i'_k\in\cx.
\tag a$$
This completes the inductive definition of the integers $c(k,\c)$.

Next we define for any $k\in[1,m]$ and any $\c\in\ce_{\ge k}$ an element $\cj_{k,\c}\in V$ by
$$\cj_{k,\c}=g(h)^\c_kg(h)^\c_{k+1}\do g(h)^\c_m\x^+$$ where 
$$g(h)^\c_\k=h(\k)^{\c(\k)}f_{i_\k}^{(\c(\k))}\text{ if }\k\in[k,m]',$$
$$g(h)^\c_\k=f_{i_\k}^{(c(\k,\c|_{\ge\k})}\text{ if } \k\in[k,m]''.$$ 
For $k\in[1,m]$ we show:
$$g(h)_kg(h)_{k+1}\do g(h)_m\x^+=\sum_{\c\in\ce_{\ge k}}\cj_{k,\c}.\tag b$$
We argue by descending induction on $k$. Assume first that $k=m$. If $k\in[1,m]'$ then  
$$g(h)_k\x^+=\sum_{n\ge0}h(k)^nf_{i_\k}^{(n)}\x^+=\sum_{\c\in\ce_{\ge k}}\cj_{k,\c},$$
as required. If $k\in[1,m]''$, then $g(h)_k\x^+=\ds_{i_k}\x^+=f_{i_k}^{(\la i_k,\l\ra)}\x^+$, see 3.3(b).

Next we assume that $k<m$ and that (b) holds for $k$ replaced by $k+1$. Let $\c'=\c_{\ge k+1}$.
By the induction hypothesis, the left hand side of (b) is equal to
$$g(h)_k\sum_{\c\in\ce_{\ge k+1}}\cj_{k+1,\c}.\tag c$$
If $k\in[1,m]'$, then clearly (c) is equal to the right hand side of (b).
If $k\in[1,m]''$, then from the induction hypothesis we see that for any $\nu\in\cx$ we have
$$(g(h)_{k+1}\do g(h)_m\x^+)_\nu=\sum_{\c\in\ce_{\ge k+1}}(\cj_{k+1,\c})_\nu
=\sum_{\c\in\ce_{\ge k+1;\nu}}\cj_{k+1,\c}$$
where $\ce_{\ge k+1;\nu}$ is the set of all $\c\in\ce_{\ge k+1}$ such that
$$\nu=\l-\sum_{\k\in[k+1,m]'}\c(\k)i'_\k-\sum_{\k\in[k+1,m]'',c(\k,\c_{\ge\k})\ge0}c(\k,\c_{\ge\k})i'_k.$$
Using this and 3.3(c) we see that
$$\align&g(h)_kg(h)_{k+1}\do g(h)_m\x^+
=\sum_{\nu\in\cx}f_{i_k}^{(\la i_k,\nu\ra)}((g(h)_{k+1}g(h)_{k+2}\do g(h)_m \x^+)_\nu)\\&=
\sum_{\nu\in\cx}f_{i_k}^{(\la i_k,\nu\ra)}\sum_{\c\in\ce_{\ge k+1;\nu}}\cj_{k+1,\c}
=\sum_{\c\in\ce_{\ge k}}f_{i_k}^{(c(k,\c))}\cj_{k+1,\c|_{\ge k+1}}=\sum_{\c\in\ce_{\ge k}}\cj_{k,\c}.
\endalign$$
This completes the inductive proof of (b).

In particular, we have
$$g(h)_1g(h)_2\do g(h)_m\x^+=\sum_{\c\in\ce}\cj_{1,\c},\tag d$$
where $\ce$ is the set of all maps $\c:[1,m]'@>>>\NN$.
This shows that for any $b\in\b$ there exists a polynomial $P_b$ in the variables $x_k,k\in[1,m]'$
with coefficients in $\NN$ such that the coefficient of $b$ in $g(h)_1g(h)_2\do g(h)_m\x^+$
is obtained by substituting in $P_b$ the variables $x_k$ by $h(k)\in\RR_{>0}$ for $k\in[1,m]',h\in A_{\ge0}$. 
Each coefficient of this polynomial is a sum of products of expressions of the form $d_{b_1,b_2,i,n}\in\NN$
(see 1.4); if one of these coefficients is $\ne0$ then after the substitution $x_k\m h(k)\in\RR_{>0}$ we obtain
an element in $\RR_{>0}$ while if all these coefficients are $0$ then the same substitution gives $0$.
Thus, there is a well defined subset $\b_{v,\ii}$ of $\b$ such that
$P_b|_{x_k=h(k)}$ is in $\RR_{>0}$ if $b\in\b_{v,\ii}$ and is $0$ if $b\in\b-\b_{v,\ii}$.

For a semifield $K_1$ we denote by $A(K_1)$ the set of maps $h:[1,m]'@>>>K_1$.
For any $h\in K_1$ we can substitute in $P_b$ the variables $x_k$ by $h(k)\in K_1$ for $k\in[1,m]'$; the result 
is an element $P_{b,h,K_1}\in K_1^!$. Clearly, we have $P_{b,h,K_1}\in K_1$ if $b\in\b_{v,\ii}$ and 
$P_{b,h,K_1}=\circ$ if $b\in\b-\b_{v,\ii}$.

From 3.2(b) we see that $b_w\in\b_{v,\ii}$. 

We see that for a semifield $K_1$, $h\m\sum_{b\in\b}P_{b,h,K_1}b$ is a map $\th_{K_1}:A(K_1)@>>>V(K_1)-\un\circ$ 
and 
$$\th_{K_1}(A(K_1))\sub\{\x\in V(K_1);\supp(\x)=\b_{v,\ii}\}.\tag d$$
($\supp(\x)$ as in 1.4.) Let $\o_{K_1}:A(K_1)@>>>P(K_1)$ be the composition of $\th_{K_1}$ with the obvious map 
$V(K_1)-\un\circ@>>>P(K_1)$. From the definitions, if $K_1@>>>K_2$ is a homomorphism of semifields, then we have 
a commutative diagram
$$\CD
A(K_1)@>\o_{K_1}>>P(K_1)\\
@VVV               @VVV\\
A(K_2)@>\o_{K_2}>>P(K_2)    
\endCD$$
where the vertical maps are induced by $K_1@>>>K_2$.

\subhead 3.5\endsubhead 
In this subsection we assume that $m\ge1$. We will consider two cases:

(I) $t_1=s_{i_1}$,

(II) $t_1=1$.
\nl
In case (I) we set $(v',w')=(s_{i_1}v,s_{i_1}w)$, $\ii'=(i_2,i_3,\do,i_m)\in\ci_{w'}$. We have $v'\le w'$
and the analogue of the sequence $q_1,q_2,\do,q_m$ in 3.2 for $(v',w',\ii')$ is $q_2,q_3,\do,q_m$.

In case (II) we set $(v',w')=(v,s_{i_1}w)$, $\ii'=\ii$. We have $v'\le w'$ and the analogue of the sequence 
$q_1,q_2,\do,q_m$ in 3.2 for $(v',w',\ii')$ is $q_2,q_3,\do,q_m$.
For a semifield $K_1$ let $A'(K_1)$ be the set of maps $[2,m]'@>>>K_1$ (notation of 3.4) and let
$\th'_{K_1}:A'(K_1)@>>>V(K_1)-\un\circ$, $\o'_{K_1}:A'(K_1)@>>>P(K_1)$ be the analogues of $\th_{K_1},\o_{K_1}$
in 3.4 when $v,w$ is replaced by $v',w'$. From the definitions, in case (I), for $h\in A(K_1)$ we have 

(a) $\th_{K_1}(h)=\ct_{i_1,K_1}(\th'_{K_1}(h|_{[2,m]'})$
\nl
(notation of 2.6(a); in this case we have $\th'_{K_1}(h|_{[2,m]'})\in V(K_1)^{e_{i_1}}$ 
by 3.3(a) and the arguments following it); hence

(b) $\o_{K_1}(h)=[\ct_{i_1,K_1}](\o'_{K_1}(h|_{[2,m]'})$
\nl
where $[\ct_{i_1,K_1}]$ is the bijection $(V(K_1)^{e_{i_1}}-\un\circ)/K_1@>>>(V(K_1)^{f_{i_1}}-\un\circ)/K_1$ 
induced by $\ct_{i_1,K_1}:V(K_1)^{e_{i_1}}@>>>V(K_1)^{f_{i_1}}$
(the image of $\o'_{K_1}(h|_{[2,m]'})$ is contained in $(V(K_1)^{e_{i_1}}-\un\circ)/K_1$).

From the definitions, in case (II), for $h\in A(K_1)$ we have 

(c) $\th_{K_1}(h)=(-i_1)^{h(i_1)}(\th'_{K_1}(h|_{[2,m]'})$
\nl
(notation of 1.4). 

\subhead 3.6\endsubhead 
In the remainder of this section we assume that $\l\in\cx^{++}$. In the setup of 3.5, let 
$h,\ti h$ be elements of $A(K_1)$. Let $\x=\th'_{K_1}(h|_{[2,m]'})$, $\ti\x=\th'_{K_1}(\ti h|_{[2,m]'})$
be such that $(-i_1)^{h(i_1)}(\x)$, $(-i_1)^{\ti h(i_1)}(\ti\x)$ have the same image in $P(K)$. We show:

(a) {\it $h(i_1)=\ti h(i_1)$ and $\x,\ti\x$ have the same image in $P(K)$.}
\nl
By 3.2(a),(b) (for $w'$ instead of $w$),

(b) {\it
$b_{w'}$ appears in $\x$ with coefficient $c\in K_1$; if $b\in\b$ appears in $\x$ with coefficient $\ne\circ$
then $\nu_b\ne\nu_{b_{w'}}+i'_1$. }

Similarly,

(c) {\it
$b_{w'}$ appears in $\ti\x$ with coefficient $\ti c\in K_1$; if $b\in\b$ appears in $\ti\x$ with coefficient 
$\ne\circ$ then $\nu_b\ne\nu_{b_{w'}}+i'_1$. }
\nl
From our assumption on $\l$ we have $b_{w'}\ne b_w=f_{i_0}^{(n)}b_{w'}$ and $f_{i_0}^{(1)}b_{w'}\ne\un\circ$. By 
(b),(c) we have

$(-i_1)^{h(i_1)}(\x)=c\b_{w'}+h(i_1)cf_{i_0}^{(1)}b_{w'}+\text{$K_1^!$-comb. of $b\in\b$ of other weights}$,

$(-i_1)^{\ti h(i_1)}(\ti\x)=\ti c\b_{w'}+\ti c\ti h(i_1)f_{i_0}^{(1)}b_{w'}
+\text{$K_1^!$-comb. of $b\in\b$ of other weights}$.
\nl
We deduce that for some $k\in K_1$ we have $\ti c=kc$, $\ti c\ti h(i_1)=kch(i_1)$. It follows that 
$h(i_1)=\ti h(i_1)$. Using this and our assumption, we see that for some $k\in K_1$ we have
$(-i_1)^{h(i_1)}(\x)=(-i_1)^{h(i_1)}(c\ti\x)$. Using 1.4(a) we deduce $\x=c\ti\x$. This proves (a).

\subhead 3.7\endsubhead 
In the setup of 3.4 we show:

(a) {\it $\o_{K_1}:A(K_1)@>>>P(K_1)$ is injective.}
\nl
We argue by induction on $m$. If $m=0$ there is nothing to prove. We now assume that $m\ge1$.
Let $\o'_{K_1}:A'(K_1)@>>>P(K_1)$ be as in 3.5. By the induction hypothesis, $\o'_{K_1}$ is injective.
In case I (in 3.5), we use 3.5(b) and the bijectivity of $[\ct_{i_1,K_1}]$ to deduce that $\o_{K_1}$ is injective.
In case II (in 3.5), we use 3.5(c) and 3.6(a) to deduce that $\o_{K_1}$ is injective. This proves (a).

\subhead 3.8\endsubhead 
According to \cite{MR}, 

(a) {\it
$h\m\s(h)B^+\s(h)\i$ defines an isomorphism $\t$ from $A$ to an open subvariety of $\cb_{v,w}$ containing 
$(\cb_{v,w})_{\ge0}$ and $\t$ restricts to a bijection $A_{\ge0}@>\si>>(\cb_{v,w})_{\ge0}$. }
\nl
(The existence of a homeomorphism $\RR_{>0}^{|w|-|v|}@>\si>>(\cb_{v,w})_{\ge0}$ was conjectured in \cite{L94b}.)

We define $\tA_{\ge0}$ in terms $A$ and its subset $A_{\ge0}$ as in 1.9. Note that $\tA_{\ge0}$ can be identified
 with the set of maps $h:[1,m]'@>>>K$ that is, with $A(K)$ (notation of 3.4). Now $\t:A@>>>\cb_{v,w}$ (see (a)) 
carries $A_{\ge0}$ onto the subset $(\cb_{v,w})_{\ge0}$ of $\cb_{v,w}$ hence it induces a map 

(b) $A(K)=\tA_{\ge0}@>>>\wt{\cb_{v,w}}_{\ge0}$ which is a bijection.
\nl
(We use (a) and 1.9(a)). 

\subhead 3.9\endsubhead 
From the definition we deduce that we have canonically

(a) $\tcb_{\ge0}=\sqc_{v,w \text{ in } W,v\le w}\wt{\cb_{v,w}}_{\ge0}$.
\nl
The left hand side is identified in 1.10 with $P^\bul(K)$, a subspace of $P(K)$.
Hence the subset $\wt{\cb_{v,w}}_{\ge0}$ of $\tcb_{\ge0}$ can be viewed as a subset $P_{v,w}(K)$
of $P(K)$ and 3.8(b) defines a bijection of $A(K)$ onto $P_{v,w}(K)$.
The composition of this bijection with the imbedding $P_{v,w}(K)\sub P(K)$ coincides with the map
$\o_K:A@>>>P(K)$ in 3.4. (This follows from definitions.)

Similarly, the composition of the imbeddings 
$$(\cb_{v,w})_{\ge0}\sub\cb_{\ge0}=P^\bul_{\ge0}\sub P_{\ge0}=P(\RR_{>0})$$
(see 1.7(a)) can be identified via 3.8(a) with the imbedding 
$\o_{\RR_{>0}}:A_{\ge0}@>>>P(\RR_{>0})$ whose image is denoted by $P_{v,w}(\RR_{>0})$.

Recall that $P^\bul(\ZZ)$ is the image of $P^\bul(K)$ under the map $P(K)@>>>P(\ZZ)$ induced by
$r:K@>>>\ZZ$ (see 1.11). For $v\le w$ in $W$ let $P_{v,w}(\ZZ)$ be the image of $P_{v,w}(K)$ under the
map $P(K)@>>>P(\ZZ)$. We have clearly $P^\bul(\ZZ)=\cup_{v\le w}P_{v,w}(\ZZ)$.
From the commutative diagram in 3.4 attached to $r:K@>>>\ZZ$ we deduce a commutative diagram
$$\CD
A(K)@>>>P_{v,w}(K)\\
@VVV               @VVV\\
A(\ZZ)@>>>P_{v,w}(\ZZ)    
\endCD$$
in which the vertical maps are surjective and the upper horizontal map is a bijection. It follows that the
lower horizontal map is surjective; but it is also injective (see 3.7(a)) hence bijective.

\subhead 3.10\endsubhead 
We return to the setup of 3.4.
If $K_1$ is one of the semifields $\RR_{>0},K,\ZZ$, then the elements of $P_{v,w}(K_1)$ are
 represented by elements of 
$\x\in V(K_1)-\un\circ$ with $\supp(\x)=\b_{v,\ii}$. In the case where $K_1=\RR_{>0}$, $P_{v,w}(K_1)$
depends only on $v,w$ and not on $\ii$. It follows that $\b_{v,\ii}$ depends only on $v,w$ not on $\ii$
hence we can write $\b_{v,w}$ insead of $\b_{v,\ii}$.

Note that in \cite{L19b, 2.4} it was conjectured that the set $[[v,w]]$ defined in \cite{L19b, 2.3(a)} in type 
$A_2$ should make sense in general. This conjecture is now established by taking $[[v,w]]=\b_{v,w}$.

Using 2.4(a) and the definitions we see that
$$\b_{v,w}\sub\b^w\cap\ph(\b^{vw_I}).\tag a$$
We expect that this is an equality (a variant of a conjecture in \cite{L19b, 2.4}, see also \cite{L19b, 2.3(a)}).
From 3.4 we see that
$$b_w\in\b_{v,w}.\tag b$$
From 2.3(d) we deduce:
$$\ph(\b_{ww_I,vw_I})=\b_{v,w}.\tag c$$
Using (b),(c) we deduce:
$$\ph(b_{vw_I})\in\b_{v,w}.\tag d$$

\subhead 3.11\endsubhead 
For $K_1$ as in 3.10 and for $v\le w$ in $W$, $v'\le w'$ in $W$, we show:

(a) {\it If $P_{v,w}(K_1)\cap P_{v',w'}(K_1)\ne\emp$, then $v=v'$, $w=w'$.}
\nl
If $K_1$ is $\RR_{>0}$ or $K$ this is already known. We will give a proof of (a) which applies also
when $K_1=\ZZ$. From the results in 3.10 we see that it is enough to show:

(b) {\it If $\b_{v,w}=\b_{v',w'}$, then $v=v'$, $w=w'$.}
\nl
From 3.10(b) we have $b_{w'}\in\b_{v',w'}$ hence $b_{w'}\in\b_{v,w}$ so that (using 3.10(a)) we have
$b_{w'}\in\b^w$. Using 2.1(a) we deduce that $b_{w'}\in V'{}^\ii$ (with $\ii$ as in 2.1).
It follows that either $b_{w'}=b_w$ or $\nu_{b_{w'}}-\nu_{b_w}$ is of the form $j'_1+j'_2+\do+j'_k$
with $j_t\in I$ and $k\ge1$. Interchanging the roles of $w,w'$ we see that
either $b_w=b_{w'}$ or $\nu_{b_w}-\nu_{b_{w'}}$ is of the form $\ti j'_1+\ti j'_2+\do+\ti j'_{k'}$
with $\ti j_t\in I$ and $k'\ge1$. If $b_w\ne b_{w'}$ then we must have 
$j'_1+j'_2+\do+j'_k+\ti j'_1+\ti j'_2+\do+\ti j'_{k'}=0$, which is absurd. Thus we have $b_w=b_{w'}$.
Since $\l\in\cx^{++}$ this implies $w=w'$.

Now applying $\ph$ to the first equality in (a) and using 3.10(c) we see that
$\b_{ww_I,vw_I}=\b_{w'w_I,v'w_I}$. Using the first part of the argument with
$v,w,v',w'$ replaced by $ww_I,vw_I,w'w_I,v'w_I$, we see that $vw_I=v'w_I$ hence $v=v'$. This completes the
proof of (b) hence that of (a).

Now the proof of Theorem 0.2 is complete.

\subhead 3.12\endsubhead 
Now $\ph:\cb@>>>\cb$ (see 2.3) induces an involution $\tcb@>>>\tcb$ and an involution $\tcb_{\ge0}@>>>\tcb_{\ge0}$
denoted again by $\ph$. From 2.3(a),(d) we deduce that this involution restricts to a bijection
$\wt{\cb_{ww_I,vw_I}}_{\ge0}@>>>\wt{\cb_{v,w}}_{\ge0}$ for any $v\le w$ in $W$. The involution
$\ph:\tcb_{\ge0}@>>>\tcb_{\ge0}$ can be viewed as an involution of $P^\bul(K)$ which coincides with the
restriction of the involution $\ph:P(K)@>>>P(K)$ in 2.7.
The last involution is compatible with the involution $\ph:P(\ZZ)@>>>P(\ZZ)$ in 2.7 under the map $P(K)@>>>P(\ZZ)$
induced by $r:K@>>>\ZZ$. It follows the image $P^\bul(\ZZ)$ of $P^\bul(K)$ under $P(K)@>>>P(\ZZ)$ is stable under
$\ph:P(\ZZ)@>>>P(\ZZ)$. Thus there is an induced involution $\ph$ on $\cb(\ZZ)=P^\bul(\ZZ)$ which carries
$P_{ww_I,vw_I}(\ZZ)$ onto $P_{v,w}(\ZZ)$ for any $v\le w$ in $W$.

\head 4. Independence on $\l$\endhead
\subhead 4.1\endsubhead
For $\l,\l'$ in $\cx^+$ let ${}^{\l,\l'}P$ be the set of lines in ${}^\l V\ot {}^{\l'}V$. We
define a linear map $E:{}^\l V\T{}^{\l'}V@>>>{}^\l V\ot{}^{\l'}V$ by $(\x,\x')\m\x\ot\x'$.
This induces a map $\bE:{}^\l P\T{}^{\l'}P@>>>{}^{\l,\l'}P$.

Let $K_1$ be a semifield. Let $\cs={}^\l\b\T{}^{\l'}\b$. 
Let ${}^{\l,\l'}V(K_1)$ be the set of formal sums $u=\sum_{s\in\cs}u_ss$ where $u_s\in K_1^!$. This is a 
monoid under addition (component by component) and we define scalar multiplication 
$$K_1^!\T{}^{\l,\l'}V(K_1)@>>>{}^{\l,\l'}V(K_1)$$
by $(k,\sum_{s\in\cs}u_ss)\m\sum_{s\in\cs}(ku_s)s$. 
Let $\End({}^{\l,\l'}V(K_1))$ be the set of maps 
$\z:{}^{\l,\l'}V(K_1)@>>>{}^{\l,\l'}V(K_1)$ such that $\z(\x+\x')=\z(\x)+\z(\x')$ for 
$\x,\x'$ in ${}^{\l,\l'}V(K_1)$ and $\z(k\x)=k\z(\x)$ for $\x\in{}^{\l,\l'}V(K_1),k\in K_1^!$. This is a monoid 
under composition of maps.

We define a map
$$E(K_1):{}^\l V(K_1)\T {}^{\l'}V(K_1)@>>>{}^{\l,\l'}V(K_1)$$
by 
$$(\sum_{b_1\in{}^\l\b}\x_{b_1}),(\sum_{b'_1\in{}^{\l'}\b}\x'_{b'_1})\m 
\sum_{(b_1,b'_1)\in\cs}\x_{b_1}\x'_{b'_1}(b_1,b'_1).$$
We define a map 
$$\End({}^\l V(K_1))\T\End({}^{\l'}V(K_1))@>>>\End({}^{\l,\l'}V(K_1))$$
 by
$(\t,\t')\m[(b_1,b'_1)\m E(K_1)(\t(b_1),\t'(b'_1))$. Composing this map with the map
$$\fG(K_1)@>>>\End({}^\l V(K_1))\T\End({}^{\l'}V(K_1))$$
whose components are the maps 
$$\fG(K_1)@>>>\End({}^\l V(K_1)),\qua\fG(K_1)@>>>\End({}^{\l'}V(K_1))$$
 in 1.5 we obtain a map $\fG(K_1)@>>>\End({}^{\l,\l'}V(K_1))$ which 
is a monoid homomorphism. Thus $\fG(K_1)$ acts on ${}^{\l,\l'}V(K_1)$; it also acts on
${}^\l V(K_1)\T{}^{\l'}V(K_1)$ (by 1.5) and the two actions are compatible with $E(K_1)$.

Let $\un\cir$ be the element $u\in{}^{\l,\l'}V(K_1)$ such that
$u_s=\circ$ for all $s\in\cs$. Let ${}^{\l,\l'}P(K_1)$ be the set of orbits of the free $K_1$ action (scalar
multiplication) on ${}^{\l,\l'}V(K_1)-\un\circ$. Now $E(K_1)$ restricts to a map
$$({}^\l V(K_1))-\un\circ)\T({}^{\l'}V(K_1)-\un\circ)@>>>{}^{\l,\l'}V(K_1)-\un\circ$$
and induces an (injective) map
$$\bE(K_1):{}^\l P(K_1)\T{}^{\l'}P(K_1)@>>>{}^{\l,\l'}P(K_1).$$
Now $\fG(K_1)$ acts naturally on ${}^\l P(K_1)\T{}^{\l'}P(K_1)$ and on ${}^{\l,\l'}P(K_1)$; these
$\fG(K_1)$-actions are compatible with $\bE(K_1)$.

\subhead 4.2\endsubhead
For $\l,\l'$ in $\cx^+$ there is a unique linear map 
$$\G:{}^{\l+\l'}V@>>>{}^\l V\ot{}^{\l'}V$$
which is compatible with the $G$-actions and takes ${}^{\l+\l'}\x^+$ to ${}^\l\x^+\ot{}^{\l'}\x^+$.
This induces a map $\bar\G:{}^{\l+\l'}P@>>>{}^{\l,\l'}P$.

For $b\in{}^{\l+\l'}\b$ we have 
$$\G(b)=\sum_{(b_1,b'_1)\in\cs}e_{b,b_1,b'_1}b_1\ot b'_1$$
where $e_{b,b_1,b'_1}\in\NN$. (This can be deduced from the positivity property \cite{L93, 14.4.13(b)} of the
homomorphism $r$ in \cite{L93, 1.2.12}.) There is a unique map 
$$\G(K_1):{}^{\l+\l'}V(K_1)@>>>{}^{\l,\l'}V(K_1)$$
compatible with addition and scalar multiplication and such that for $b\in{}^{\l+\l'}\b$ we have
$$\G(K_1)(b)=\sum_{(b_1,b'_1)\in\cs}e_{b,b_1,b'_1}(b_1,b'_1)$$
where $e_{b,b_1,b'_1}$ are viewed as elements of $K_1^!$. Since $\G$ is injective, for any 
$b\in{}^{\l+\l'}\b$ we have $e_{b,b_1,b'_1}\in\NN-\{0\}$ for some 
$b_1,b'_1$, hence $e_{b,b_1,b'_1}\in K_1$, when viewed as an element of $K_1^!$. It follows that 
$\G(K_1)$ maps ${}^{\l+\l'}V(K_1)-\un\circ$ into ${}^{\l,\l'}V(K_1)-\un\circ$. Hence 
$\G(K_1)$ defines an (injective)  map
$$\bar\G(K_1):{}^{\l+\l'}P(K_1)@>>>{}^{\l,\l'}P(K_1)$$
which is compatible with the action of $\fG(K_1)$ on the two sides.

\subhead 4.3\endsubhead
We now assume that $K_1$ is either $K$ as in 0.1(i) or $\ZZ$ as in 0.1(ii) and that 
$\l\in\cx^{++},\l'\in\cx^+$ so that $\l+\l'\in\cx^{++}$. We have the following result.

(a) {\it Let $\cl\in{}^{\l+\l'}P^\bul(K_1)$. Then $\bar\G(K_1)(\cl)=\bE(K_1)(\cl_1,\cl'_1)$ for some
$(\cl_1,\cl'_1)\in{}^\l P^\bul(K_1)\T{}^{\l'}P(K_1)$ (which is unique, by the injectivity of $\bE(K_1)$). 
Thus, $\cl\m\cl_1$ is a well defined map $H(K_1):{}^{\l+\l'}P^\bul(K_1)@>>>{}^\l P^\bul(K_1)$.}
\nl
We shall prove (a) for $K_1=\ZZ$ assuming that it is true for $K_1=K$.
We can find $\ti\cl\in{}^{\l+\l'}P^\bul(K)$ such that $\cl\in{}^{\l+\l'}P^\bul(\ZZ)$  is the image of
$\ti\cl$ under the map ${}^{\l+\l'}P^\bul(K)@>>>{}^{\l+\l'}P^\bul(\ZZ)$ induced by $r:K@>>>\ZZ$.
By our assumption we have $\bar\G(K)(\ti\cl)=\bE(K)(\ti\cl_1,\ti\cl'_1)$ with
$(\ti\cl_1,\ti\cl'_1)\in{}^\l P^\bul(K)\T{}^{\l'}P(K)$. Let $\cl_1$ (resp. $\cl'_1$)
be the image of $\ti\cl_1$ (resp. $\ti\cl'_1$) under the map ${}^\l P^\bul(K)@>>>{}^\l P^\bul(\ZZ)$ 
(resp. ${}^{\l'}P(K)@>>>{}^{\l'}P(\ZZ)$) induced by $r:K@>>>\ZZ$. From the definitions we see that
$\bar\G(\ZZ)(\cl)=\bE(\ZZ)(\cl_1,\cl'_1)$. This proves the existence of  $(\cl_1,\cl'_1)$. The proof of (a) in 
the case where $K_1=K$ will be given in 4.6.

Assuming that (a) holds, we have a commutative diagram
$$\CD
{}^{\l+\l'}P^\bul(K)@>H(K)>>{}^\l P^\bul(K)\\
@VVV                        @VVV   \\
{}^{\l+\l'}P^\bul(\ZZ)@>H(\ZZ)>>{}^\l P^\bul(\ZZ)
\endCD$$
in which the vertical maps are induced by $r:K@>>>\ZZ$. 

\subhead 4.4\endsubhead
We preserve the setup of 4.3.
For each $w\in W$ we assume that a sequence $\ii_w=(i_1,i_2,\do,i_m)\in\ci_w$ has been chosen (here $m=|w|$).
Let $\cz(K_1)=\sqc_{v\le w\text{ in }W}A_{v,w}(K_1)$ where 
$A_{v,w}(K_1)$ is the set of all maps $[1,m]'@>>>K_1$ (with $[1,m]'$ defined as in 3.2 in terms of $v,w$ and
$\ii=\ii_w$). From the results in 3.9 we have a bijection
$${}^\l D(K_1):\cz(K_1)@>\si>>{}^\l P^\bul(K_1)$$
whose restriction to $A_{v,w}(K_1)$ is as in the last commutative diagram in 3.9 (with $\ii=\ii_w$). Replacing 
here $\l$ by $\l+\l'$ we obtain an analogous bijection
$${}^{\l+\l'}D(K_1):\cz(K_1)@>\si>>{}^{\l+\l'} P^\bul(K_1).$$
From the commutative diagram in 3.4 we deduce a commutative diagram
$$\CD
\cz(K)@>{}^\l D(K)>>{}^\l P^\bul(K)\\
@VVV            @VVV    \\
\cz(\ZZ)@>{}^\l D(\ZZ)>>{}^\l P^\bul(\ZZ)
\endCD$$
and a commutative diagram
$$\CD
\cz(K)@>{}^{\l+\l'}D(K)
>>{}^{\l+\l'} P^\bul(K)\\
@VVV            @VVV    \\
\cz(\ZZ)@>{}^{\l+\l'}D(\ZZ)>>{}^{\l+\l'} P^\bul(\ZZ)
\endCD$$
in which the vertical maps are induced by $r:K@>>>\ZZ$.

\subhead 4.5\endsubhead
We preserve the setup of 4.3. We assume that 4.3(a) holds.
From the commutative diagrams in 4.3, 4.4 we deduce a commutative diagram
$$\CD
\cz(K)@>({}^\l D(K))\i H(K){}^{\l+\l'}D(K)>>    \cz(K)\\
@VVV            @VVV    \\
\cz(\ZZ)@>({}^\l D(\ZZ))\i H(\ZZ){}^{\l+\l'}D(\ZZ)>>    \cz(\ZZ)
\endCD$$
in which the vertical maps are induced by $r:K@>>>\ZZ$.
Recall that $K_1$ is $K$ or $\ZZ$. We have the following result.

(a) {\it $({}^\l D(K_1))\i H(K_1){}^{\l+\l'}D(K_1)$ is the identity map $\cz(K_1)@>>>\cz(K_1)$.}
\nl
If (a) holds for $K_1=K$ then it also holds for $K_1=\ZZ$, in view of the commutative diagram above in
which the vertical maps are surjective. The proof of (a) in the case $K_1=K$ will be given in 4.7.

From (a) we deduce:

(b) {\it $H(K_1)$ is a bijection.}

\subhead 4.6\endsubhead
In this subsection we assume that $K_1=K$. Let $\kk=\CC(x)$ where $x$ is an indeterminate. We have $K^!\sub\kk$. 
For any $\l\in\cx^+$ we set ${}^\l V_\kk=\kk\ot{}^\l V$. This is naturally a module over the group $G(\kk)$ of 
$\kk$ points of $G$.
Let $\cb(\kk)$ be the set of subgroups of $G(\kk)$ that are $G(\kk)$-conjugate to $B^+(\kk)$, the group of
$\kk$-points of $B^+$.
We identify ${}^\l V(K)$ with the set of vectors in ${}^\l V_\kk$ whose coordinates in the $\kk$-basis
${}^\l\b$ are in $K^!$. In the case where $\l\in\cx^{++}$, we identify ${}^\l V^\bul(K)-0$ with the set of all 
$\x\in{}^\l V(K)-0$ such that the stabilizer in $G(\kk)$
of the line $[\x]$ belongs to $\cb(\kk)$. (For a nonzero vector $\x$ in a $\kk$-vector space we denote by $[\x]$
the $\kk$-line in that vector space that contains $\x$.)

Now let $\l\in\cx^{++},\l'\in\cx^+$. We show that 4.3(a) holds for $\l,\l'$. We identify ${}^{\l,\l'}V(K)$ with 
the set of vectors in ${}^\l V_\kk\ot_\kk{}^{\l'}V_\kk$ whose coordinates in the $\kk$-basis 
${}^{\l}\b\ot{}^{\l'}\b$ are in $K^!$. 

Then $E(K)$ becomes the restriction of the homomorphism of $G(\kk)$-modules
$E':{}^\l V_\kk \T {}^{\l'}V_\kk@>>>{}^\l V_\kk\ot_\kk{}^{\l'}V_\kk$ given by $(\x,\x')\m\x\ot_\kk\x'$
and $\G(K)$ becomes the restriction of the homomorphism of $G(\kk)$-modules
$\G':{}^{\l+\l'}V_\kk@>>>{}^\l V_\kk\ot_\kk{}^{\l'}V_\kk$
obtained from $\G$ by extension of scalars.

Let $L_\l=[{}^\l\x^+]\sub{}^\l V_\kk$, $L_{\l'}=[{}^{\l'}\x^+]\sub{}^{\l'}V_\kk$,
$L_{\l+\l'}=[{}^{\l+\l'}\x^+]\sub{}^{\l+\l'}V_\kk$. Now let $\x\in{}^{\l+\l'}V^\bul(K)-0$.
Then $[\x]=gL_{\l+\l'}$ for some $g\in G(\kk)$ hence
$$\G'([\x])=g(L_\l\ot L_{\l'})=(gL_\l)\ot(g(L_{\l'})=E'(gL_\l,g(L_{\l'})=E'([g({}^\l\x^+)],[g({}^{\l'}\x^+)]).$$
To prove 4.3(a) in our case it is enough to prove that for some $c,c'$ in $\kk^*$ we have
$cg({}^\l\x^+)\in{}^\l V(K)$, $c'g({}^{\l'}\x^+)\in{}^{\l'}V(K)$.
We have $\x=c_0g({}^{\l+\l'}\x^+)$ for some $c_0\in\kk^*$ and $\G'(\x)=\G(\x)\in{}^{\l,\l'}V(K)$. Thus,
 $c_0\G'(g({}^{\l+\l'}\x^)\in{}^{\l,\l'}V(K)$ that is,
$c_0 (g{}^\l\x^+)\ot(g{}^{\l'}\x^+)\in{}^{\l,\l'}V(K)$.
It is enough to show: 

(a) {\it If $z\in{}^\l V_\kk$, $z'\in {}^{\l'}V_\kk$, $c_0\in\kk^*$ satisfy $c_0z\ot z'\in{}^{\l,\l'}V(K)-0$, 
then $cz\in{}^\l V(K)-0$, $c'z'\in {}^{\l'}V(K)-0$ for some $c,c'$ in $\kk^*$.}
\nl
We write $z=\sum_{b\in{}^\l\b}z_bb$, $z'=\sum_{b'\in{}^{\l'}\b}z'_{b'}b'$ with $z_b,z'_{b'}$ in $\kk$.
We have $c_0z_bz'_{b'}\in K^!$ for all $b,b'$. Replacing $z$ by $c_0z$ we can assume that $c_0=1$ so that
$z_bz'_{b'}\in K^!$ for all $b,b'$ and $z_bz'_{b'}\ne0$ for some $b,b'$.
Thus we can find $b'_0\in{}^{\l'}\b$ such that $z'_{b'_0}\in K$. We have
$z_bz'_{b'_0}\in K^!$ for all $b$. Replacing $z$ by $z'_{b'_0}z$ we can assume that
$z_b\in K^!$ for all $b$. We can find $b_0\in{}^\l\b$ such that $z_{b_0}\in K$. We have
$z_{b_0}z'_{b'}\in K^!$ for all $b'$. It follows that
$z'_{b'}\in K^!$ for all $b'$. This proves (a) and completes the proof of 4.3(a).

\subhead 4.7\endsubhead
We preserve the setup of 4.3 and assume that $K_1=K$. We show that 4.5(a) holds in this case.
Let $v\le w,\ii$ be as in 3.2 and let $A(K_1)$ be as in 3.4. Let $h\in A(K_1)$. 
We have ${}^{\l+\l'}D(K_1)(h)=[\s_{K_1}(h){}^{\l+\l'}\x^+]$
where $\s_{K_1}:A(K_1)@>>>G(\kk)$ is defined by the same formula as $\s$ in 3.2.
(Note that for $i\in I$, $y_i(t)\in G(\kk)$ is defined for any $t\in\kk$.) Hence 
$$\align&\bar\G(K_1){}^{\l+\l'}D(K_1)(h)=[(\s_{K_1}(h){}^\l\x^+)\ot(\s_{K_1}(h){}^{\l'}\x^+)]\\&
=\bE(K_1)([\s_{K_1}(h){}^\l\x^+],[\s_{K_1}(h){}^{\l'}\x^+])\endalign$$
so that
$$H(K_1){}^{\l+\l'}D(K_1)(h)=[\s_{K_1}(h){}^\l\x^+]={}^\l D(K_1)(h).$$
This shows that the map in 4.5(a) takes $h$ to $h$ for any $h\in A(K_1)$. This proves 4.5(a).

\subhead 4.8\endsubhead
We now assume that $K_1$ is either $K$ as in 0.1(i) or $\ZZ$ as in 0.1(ii) and that 
$\l\in\cx^{++},\l'\in\cx^{++}$.
From 4.3(a),4.5(a) we have a well defined bijection
$H(K_1):{}^{\l+\l'}P^\bul(K_1)@>\si>>{}^\l P^\bul(K_1)$. Interchanging $\l,\l'$ we obtain a bijection
$H'(K_1):{}^{\l+\l'}P^\bul(K_1)@>\si>>{}^{\l'}P^\bul(K_1)$. Hence we have a bijection
$$\g_{\l,\l'}=H'(K_1)H(K_1)\i:{}^\l P^\bul(K_1)@>\si>>{}^{\l'}P^\bul(K_1).$$
From the definitions we see that $H(K_1)$ is compatible with the $\fG(K_1)$-actions. Similarly,
$H'(K_1)$ is compatible with the $\fG(K_1)$-actions. It follows that $\g_{\l,\l'}$ is compatible with the 
$\fG(K_1)$-actions. From the definitions we see that if $\l''$ is third element of $\cx^{++}$, we have 
$$\g_{\l,\l''}= \g_{\l',\l''}\g_{\l,\l'}.$$
This shows that our definition of $\cb(K_1)$ is independent of the choice of $\l$.

\head 5. The non-simplylaced case\endhead
\subhead 5.1\endsubhead
Let $\d:G@>>>G$ be an automorphism of $G$ such that $\d(B^+)=B^+,\d(B^-)=B^-$ and $\d(x_i(t))=x_{i'}(t)$,
$\d(y_i(t))=y_{i'}(t)$ for all $i\in I,t\in\CC$ where $i\m i'$ is a permutation of $I$ denoted again by 
$\d$. We define an automorphism of $W$ by $s_i\m s_{\d(i)}$ for all $i\in I$; we denote this automorphism
again by $\d$. We assume further that $s_is_{\d(i)}=s_{\d(i)}s_i$ for any $i\in I$. The fixed point set
$G^\d$ of $\d:G@>>>G$ is a connected simply connected semisimple group over $\CC$.
The fixed point set $W^\d$ of $\d:W@>>>W$ is the Weyl group of $G^\d$ and as such it has a length function
$w\m|w|_\d$.

Now $\d$ takes any Borel subgroup of $G$ to a Borel subgroup of $G$ hence it defines an automorphism of $\cb$
denoted by $\d$, with fixed point set denoted by $\cb^\d$. This automorphism restricts to a
bijection $\cb_{\ge0}@>>>\cb_{\ge0}$. We can identify $\cb^\d$ with the
 flag manifold of $G^\d$ by $B\m B\cap G^\d$. Under this identification, the totally positive part of the flag 
manifold of $G^\d$ (defined in \cite{L94b}) becomes $\cb^\d_{\ge0}=\cb_{\ge0}\cap\cb^\d$.
For $\l\in\cx$ we define $\d(\l)\in\cx$ by $\la\d(i),\d(\l)\ra=\la i,\l\ra$ for all $i\in I$.
In the setup of 1.4 assume that $\l\in\cx^{++}$ satisfies $\d(\l)=\l$. There is a unique linear isomorphism
$\d:V@>>>V$ such that $\d(g\x)=\d(g)\d(\x)$ for any $g\in G,\x\in V$ and such that $\d(\x^+)=\x^+$.
This restricts to a bijection $\b@>>>\b$ denoted again by $\d$. For any semifield $K_1$ we define a bijection
$V(K_1)@>>>V(K_1)$ by $\sum_{b\in\b}\x_bb\m\sum_{b\in\b}\x_{\d\i(b)}b$ where $\x_b\in K_1^!$. This induces a 
bijection $P(K_1)@>>>P(K_1)$ denoted by $\d$.
We now assume that $K_1$ is as in 0.1(i),(ii). Then the subset $P^\bul(K_1)$ of $P(K_1)$ is defined and is 
stable under $\d$; let $P^\bul(K_1)^\d$ be the fixed point set of $\d:P^\bul(K_1)@>>>P^\bul(K_1)$. 
Recall that $\fG(K_1)$ acts naturally on $P(K_1)$. This restricts to an action on $P^\bul(K_1)^\d$ of the monoid 
$\fG(K_1)^\d$ (the fixed point set of the isomorphism $\fG(K_1)@>>>\fG(K_1)$ induced by $\d$) which is the same
as the monoid associated in \cite{L19a} to $G^\d$ and $K_1$. We set $\cb^\d(K_1)=P^\bul(K_1)^\d$.

The following generalization of Theorem 0.2 can be deduced from Theorem 0.2.

(a) {\it The set $\cb^\d(\ZZ)$ has a canonical partition into pieces $P_{v,w;\d}(\ZZ)$ indexed by the
pairs $v\le w$ in $W^\d$. Each such piece $P_{v,w;\d}(\ZZ)$ is in bijection with $\ZZ^{|w|_\d-|v|_\d}$; in fact, 
there is an explicit bijection $\ZZ^{|w|_\d-|v|_\d}@>\si>>P_{v,w;\d}(\ZZ)$ for any reduced expression of $w$
in $W^\d$.}


\widestnumber\key{L90b}
\Refs
\ref\key{L90a}\by G.Lusztig\paper Canonical bases arising from quantized enveloping algebras
\jour J. Amer. Math. Soc.\vol3\yr1990\pages 447-498\endref
\ref\key{L90b}\by G.Lusztig\paper Canonical bases arising from quantized enveloping algebras,II
\jour Progr.of Theor. Phys. Suppl.\vol102\yr1990\pages 175-201\endref
\ref\key{L93}\by G.Lusztig\book Introduction to quantum groups\bookinfo Progr.in Math.110\publ 
Birkh\"auser Boston\yr 1993\endref
\ref\key{L94a}\by G.Lusztig\paper Problems on canonical bases\inbook Algebraic groups and their 
generalizations: quantum and infinite dimensionalmethods\bookinfo Proc. Symp. Pure Math. 56(2)
\publ Amer. Math. Soc.\yr1994\pages 169-176\endref
\ref\key{L94b}\by G.Lusztig\paper Total positivity in reductive groups\inbook Lie theory and geometry\bookinfo
 Progr.in Math. 123\publ Birkh\"auser Boston \yr1994\pages 531-568\endref
\ref\key{L98}\by G.Lusztig\paper Total positivity in partial flag manifolds\jour Represent.Th.\vol2\yr1998\pages
70-78\endref
\ref\key{L18}\by G.Lusztig\paper Positive structures in Lie theory\jour arxiv:1812.09313\toappear\endref
\ref\key{L19a}\by G.Lusztig\paper Total positivity in reductive groups,II \jour Bull. Inst. Math. Acad. Sin.
(N.S.)\vol14\yr2019\pages403-460\endref
\ref\key{L19b}\by G.Lusztig\paper Total positivity in Springer fibres\jour arxiv:1909.00896\endref
\ref\key{MR}\by R.Marsh and K.Rietsch\paper Parametrizations of flag varieties\jour Represent.Th.\vol8\yr2004
\pages212-242\endref
\endRefs
\enddocument